\documentclass[11pt]{amsart}
\usepackage{amssymb}

\newtheorem{theorem}{Theorem}[section]
\newtheorem{lemma}[theorem]{Lemma}

\newtheorem{definition}[theorem]{Definition}

\newtheorem{corollary}[theorem]{Corollary}
\newtheorem{proposition}[theorem]{Proposition}

\newtheorem{lem-def}[theorem]{Lemma-Definition}

\DeclareRobustCommand\longtwoheadrightarrow
{\relbar\joinrel\twoheadrightarrow}

\renewenvironment{proof}{{\bfseries Proof.}}{\qed}
\newcommand{\Ha}{\mathbb H}

\newcommand{\N}{\mathbb N}
\newcommand{\Z}{\mathbb Z}
\newcommand{\Q}{\mathbb Q}

\newcommand{\F}{\mathbb F}

\def\op{\operatorname}

\def\aa{\mathcal{A}}

\def\al{\alpha}

\def\ars#1{\renewcommand\arraystretch{#1}}

\def\bs{\vskip.5cm}
\def\be{\beta}

\def\chr{\op{char}}

\def\defn{\nn{\bf Definition. }}
\def\dep{\op{\mbox{\rm\small depth}}}

\def\dgk{\deg_K}

\def\diso{\lower.4ex\hbox{$\downarrow$}\raise.4ex\hbox{\mbox{\scriptsize
$\wr$}}}

\def\dta{\delta}

\def\e{\medskip}

\def\ep#1{\exp(\Pi i#1)}
\def\ep{\epsilon}

\def\fsep{f_{\op{sep}}}

\def\g{\Gamma}
\def\ga{\gamma}

\def\gal{\op{Gal}}

\def\gg{\mathcal{G}}

\def\ggm{\mathcal{G}_\mu}

\def\gi{\g_{\infty}}

\def\gm{\g_\mu}

\def\ii{i\in I}

\def\imp{\ \Longrightarrow\ }

\def\inm{\op{in}_\mu}
\def\inn{\op{in}}

\def\irr{\op{Irr}}
\def\ism{\lower.3ex\hbox{\ars{.08}$\begin{array}{c}\,\to\\\mbox{\tiny $\sim\,$}\end{array}$}}
\def\iso{\ \lower.3ex\hbox{\ars{.08}$\begin{array}{c}\lra\\\mbox{\tiny $\sim\,$}\end{array}$}\ }

\def\ka{\kappa}

\def\kb{\overline{K}}

\def\km{k_\mu}

\def\kp{\op{KP}}

\def\kpi{\op{KP}_\infty}
\def\kpm{\op{KP}(\mu)}

\def\kx{K[x]}

\def\la{\lambda}
\def\La{\Lambda}

\def\lg{l\raise.6ex\hbox to.2em{\hss.\hss}l}

\def\lra{\,\longrightarrow\,}

\def\md#1{\; \mbox{\rm(mod }{#1})}

\def\mlv{Mac Lane--Vaqui\'e\ }
\def\mm{{\mathcal M}}

\def\mub{\bar{\mu}}

\def\nn{\noindent}

\def\om{\omega}

\def\orb{\hbox to  .3em{$\backslash$}\backslash}
\def\ord{\op{ord}}
\def\p{\mathfrak{p}}

\def\phsep{\phi_{\op{sep}}}

\def\ram{^{\mbox{\tiny ram}}}

\def\ri{\rho_i}

\def\rj{\rho_j}

\def\sep{\op{sep}}
\def\sg{\sigma}

\def\sii{\ \Longleftrightarrow\ }

\def\sub{\subseteq}
\def\supp{\op{supp}}

\def\t{\theta}
\def\tame{^{\mbox{\tiny tame}}}

\def\ttt{\mathcal{T}}
\def\ty{\mathbf{t}}

\def\z{\op{Z}}

\newcounter{cs}
\stepcounter{cs}
\newcommand{\casos}{\begin{itemize}}
\newcommand{\fcasos}{\end{itemize}\setcounter{cs}{1}}

\newfont{\tit}{cmr12 scaled \magstep3}

\title{Okutsu sequences in Henselian valued fields}%with defect

\makeatletter
\@namedef{subjclassname@2010}{%
  \textup{2010} Mathematics Subject Classification}
%, 12J10}

\author[Nart]{Enric Nart}
\address{Departament de Matem\`{a}tiques,         Universitat Aut\`{o}noma de Barcelona,         Edifici C, E-08193 Bellaterra, Barcelona, Catalonia}
\email{enric.nart@uab.cat}

\thanks{Partially supported by grant PID2020-116542GB-I00  funded by the Spanish MCIN/AEI}


\begin{document}

\begin{abstract}
For $(K,v)$ a Henselian  valued field, let $\t\in\kb$ with minimal polynomial $F$ over $K$.  Okutsu sequences of $\t$ have been defined only when the extension $K(\t)/K$ is defectless. In this paper, we extend this concept to arbitrary $\t\in\kb$ and we show that these objects are essentially equivalent  to Okutsu frames of $F$ and to Mac Lane-Vaqui\'e chains of the natural valuation on $\kx$ induced by $\t$.
\end{abstract}

\maketitle

%\begin{center}\sl Preliminary version \end{center}

%\tableofcontents

\section*{Introduction}

Let $(K,v)$ be a Henselian valued field and let us still denote by $v$ the unique extension of $v$ to some fixed algebraic closure $\kb$. Let $vK=v(K^*)$ be the value group and $K\!v$ the residue field of $v$ over $K$. Note that $\g:=v\kb$ is a divisible hull of $vK$.

Consider  a finite simple field extension  $L/K$. The ramification index $e(L/K)$, inertia degree $f(L/K)$ and defect $d(L/K)$ of this extension are three natural numbers linked by 
\[
 [L\colon K]=e(L/K)f(L/K)d(L/K).
\]
By a celebrated resut of Ostrowski, $d(L/K)$ is a power of  the residual characteristic $p$ of $(K,v)$; that is, $p=1$ if $K\!v$ has characteristic zero, and $p=\chr(K\!v)$ otherwise.

Let  $L=K(\t)$ for some $\t\in\kb$ with minimal polynomial $F$ over $K$. Such an $F$ is said to be a ``generator" of $L/K$. Our leitmotif is the following relevant question.\e

\nn{\bf Problem. }\emph{Compute $e(L/K)$, $f(L/K)$, $d(L/K)$, in terms of any given generator $F$ of $L/K$.}\e

This problem was  solved by Mac Lane when $v$ is discrete of rank one and $L/K$ is separable \cite{mcla, mclb}.  
It is well-known that $d(L/K)=1$ in this 
classical situation. Seventy years later, Vaqui\'e extended Mac Lane's theory to valued fields $(K,v)$ of arbitrary rank \cite{Vaq}. 
%Let us briefly sketch this Mac Lane-Vaqui\'e theory.

%\subsection*{Mac Lane-Vaqui\'e chains of valuations on the polynomial ring}

MacLane realized that the link between $(L/K,v)$ and $F$ could be described in terms of the following valuation on the polynomial ring:
\[
v_F\colon \kx\lra \g\cup\{\infty\},\qquad g\longmapsto v_F(g)=v(g(\t)).
\]  
%By the Henselian property, $v_F$ does not depend on the choice of $\t$ among the roots of $F$ in $\kb$.  
Since $v_F^{-1}(\infty)=F\kx$, the valuations  $v_F$ and $v$ (on $L$) determine one each other through
\[
v_F\colon \kx \longtwoheadrightarrow \kx/(F)\stackrel{\sim}\lra L\stackrel{v}\lra \g\cup\{\infty\}, 
\]
where the isomorphism $\kx/(F)\stackrel{\sim}\to L $ is induced by $x\mapsto \t$.

The set of all extensions of $v$ to $\kx$ taking values inside $\g$ is partially ordered, with respect to the following ordering:
\begin{equation}\label{ordering}	
\mu\le \nu\ \sii\ \mu(f)\le\nu(f)\quad\mbox{for all }\ f\in \kx. 
\end{equation}

MacLane showed that $v_F$ can be constructed from $v$ by applying a finite number of \textbf{augmentations} of valuations on $\kx$ extending $v$:
\begin{equation}\label{mlv}
v\ \stackrel{\phi_0,\ga_0}\lra\ 	\mu_0\ \stackrel{\phi_1,\ga_1}\lra\  \mu_1\ \lra\ \cdots
	\ \stackrel{\phi_{r-1},\ga_{r-1}}\lra\ \mu_{r-1}
	\ \stackrel{F,\infty}\lra\ \mu_{r}=v_F,
\end{equation}
where $\phi_0,\dots,\phi_{r-1}\in\kx$ are certain \textbf{key polynomials} and $\ga_0,\dots,\ga_{r-1}\in vL\sub\g$. The initial step $v\to\mu_0$ is just  a formal augmentation. The monic polynomial $\phi_0$ has degree one  and $\mu_0$ is  determined by $\mu_0(\phi_0)=\ga_0$ and
\begin{equation}\label{mu0}
\mu_0\left(\sum\nolimits_{0\le n}a_n\phi_0^n\right):=\min_{0\le n}\{v(a_n)+n\ga_0\}=\min_{0\le n}\{\mu_0\left(a_n\phi_0^n\right)\}.
\end{equation}

Similarly, given a valuation $\mu$ on $\kx$, a key polynomial $\phi$ for $\mu$ (cf. Definition \ref{keypolyno}) and a value $\ga\in \g$ such that $\ga>\mu(\phi)$, then the augmentation \[\mu\stackrel{\phi,\ga}\lra\nu\] is defined as follows in terms of $\phi$-expansions. Every $g\in\kx$ can be written in a unique way as  $g=\sum\nolimits_{0\le n}a_n\phi^n$, for some $a_n\in\kx$ with $\deg(a_n)<\deg(\phi)$. Then, we define 
\begin{equation}\label{ordinary}
\nu(g):=\min_{0\le n}\{\mu(a_n)+n\ga\}=\min_{0\le n}\{\nu\left(a_n\phi^n\right)\}.
\end{equation}
In each augmentation $\mu\to\nu$, we consider some relative ramification index $e(\mu\to\nu)$ and inertia degree $f(\mu\to\nu)$. Under certain technical condition on the chain (\ref{mlv}), one has
\begin{equation}\label{ef}
\begin{array}{ccl}
e(L/K)&=&e(v\to\mu_0)\cdots e(\mu_{r-2}\to\mu_{r-1}),\\
f(L/K)&=&f(\mu_0\to\mu_1)\cdots f(\mu_{r-1}\to\mu_r).	
\end{array}	
\end{equation}

The paper \cite{OM} presents an algorithm computing the chain (\ref{mlv}) and the relative ramification index and inertia degree of each augmentation.\e

This approach of Mac Lane was generalized by Vaqui\'e to valuations of arbitrary rank. In this general setting, it is necessary to consider a new type of \textbf{limit augmentations}.

A \textbf{Mac Lane-Vaqui\'e (MLV) chain} of $v_F$ is a chain as in (\ref{mlv}) made of a mixture of ordinary augmentations (defined as in (\ref{ordinary})) and limit augmentations. A certain technical condition is imposed too, in order to ensure the validity of (\ref{ef}) and guarantee a certain unicity property. The precise definitions can be found in Section \ref{subsecMLV}.

Besides proving the existence of MLV chains, Vaqui\'e showed that the defectless extensions are precisely those where no limit augmentations appear \cite{Vaq2}. More precisely, the defect can be expressed as well as a product of relative defects:  
\[
d(L/K)=d(\mu_0\to\mu_1)\cdots d(\mu_{r-1}\to\mu_r), 
\]
where $d(\mu_n\to\mu_{n+1})=1$ if and only if the augmentation  $\mu_n\to\mu_{n+1}$ is ordinary. This ``defect formula" of Vaqui\'e was extended in \cite{NN} to arbitrary (not necessarily Henselian) valued fields $(K,v)$.

Let us emphasize that it is an open problem to design an algorithm computing MLV chains of a given $v_F$ and the relative data $e(\mu_n\to\mu_{n+1})$, $f(\mu_n\to\mu_{n+1})$, $d(\mu_n\to\mu_{n+1})$ of each augmentation. In \cite{OM}, some algorithms are described solving  this question in a few particular cases that go beyond the classical discrete rank-one setting. \e

%\subsection*{Okutsu frames and sequences}\label{subsecOKfrseq}

\textbf{Okutsu sequences} of algebraic elements were introduced by Okutsu for complete, discrete, rank-one valued fields, as a tool to construct integral bases \cite{Ok}.
An Okutsu sequence of a given $\t\in\kb$ is a finite list of algebraic elements with increasing degree over $K$:
\[
 [\al_0,\al_1,\dots,\al_{r-1},\al_r=\t],\qquad 1=\deg_K\al_0<\cdots<\deg_K\al_{r-1}<\deg_K\t,
\]
satisfying the following conditions for all $\be\in\kb$ and all $0\le i<r$:
\begin{itemize}
 \item \ $\deg_K\be<\deg_K\al_{i+1}\ \imp\ v(\t-\be)\le v(\t-\al_i)$,
 \item \ $\deg_K\be<\deg_K\al_i\ \imp\ v(\t-\be)< v(\t-\al_i)$.
\end{itemize}
Roughly speaking, the $\al_i$ are a kind of ``better approximations" to $\t$ in the $v$-adic topology, according to their degree over $K$. 

Under the form of \textbf{distinguished  chains} of algebraic elements, these objects have been developed by several authors, mainly in the Henselian case \cite{PZ, AK1, AK2}. The equivalence between Okutsu sequences and distinguished chains was established  in \cite{OkS} for defectless extensions of Henselian valued fields of an arbitrary rank. For defect extensions, neither Okutsu sequences nor distinguished chains exist.\e

Let $\irr(K)$ be the set of all monic, irreducible polynomials in $\kx$. 
The connection between Okutsu's and Mac Lane's ideas was first detected in \cite{okutsu}. In the discrete rank-one case, it was proved that \textbf{Okutsu frames} of a monic  $F\in\irr(K)$ are essentially the same objects as MLV chains of $v_F$.
Consider the following \textbf{weight} function with respect to $F$:
\[%\begin{equation}\label{weight}
 w(g):=v_F(g)/\deg(g)\in\g,\quad\mbox{for all }\ g\in\kx\setminus K.
\]%\end{equation}
An Okutsu frame of $F$ is a finite list of polynomials in $\irr(K)$, with increasing degree:
\[
 [F_0,F_1,\dots,F_{r-1},F_r=F],\qquad 1=\deg(F_0)<\cdots<\deg(F_{r-1})<\deg(F),
\]
satisfying the following conditions for all monic $g\in\kx$ and all $0\le i<r$:
\begin{itemize}
 \item \ $\deg(g)<\deg(F_{i+1})\ \imp\ w(g)\le w(F_i)$,
 \item \ $\deg(g)<\deg(F_i)\ \imp\ w(g)< w(F_i)$.
\end{itemize}

The main result of \cite{okutsu} states that the polynomials in an Okutsu frame of $F$ are key polynomials of a certain MLV chain of $v_F$, and viceversa. This result was generalized in \cite{defless} to defectless extensions of valued fields of arbitrary rank. 

Finally, in \cite{OkF} this equivalence has been extended to 
arbitrary extensions of valued fields of arbitrary rank. For defect extensions, classical Okutsu frames defined as above do not exist; thus, the definition of an Okutsu frame needs to be modified. For certain degrees, individual polynomials must be replaced with infinite families of polynomials.  \e

Now, let us sketch the content of the paper.
In Section \ref{secMLVOF}, we introduce the necessary background on MLV chains and their connection with Okutsu frames. Most of the content of this section is taken from \cite{MLV} and \cite{OkF}.

In Section \ref{secWD}, we discuss the relationship between weights and distances of polynomials in $\kx$, with respect to a fixed $\t\in\kb$. This material is crucial to prove our main result in Section \ref{secOS}, where we introduce Okutsu sequences of an arbitrary $\t\in\kb$, not necessarily defectless. In Theorem \ref{SF}, we show that they are essentially equivalent to Okutsu frames of the minimal polynomial of $\t$ over $K$.

In Section \ref{secTame}, we review (and correct) the results of \cite{OkS} on the comparison of the \textbf{main invariant} $\dta(\t)$ and  \textbf{Krasner's constant} $\om(\t)$, defined as follows
\[\ars{1.2}
\begin{array}{l}
\dta(\t):=\dta(F)=\sup\{v(\t-\al)\mid \al\in\kb,\ \deg_K\al<\deg_K\t\},\\
\om(\t):=\om(F)=\max\{v(\t-\t')\mid \t'\in \z(F),\ \t'\ne\t\},
\end{array}
\]
where $\z(F)$ is the multiset of all roots of $F$ in $\kb$, counting multiplicities.
By Krasner's lemma, $\dta(\t)\le\om(\t)$, if $\t$ is separable. It is well-known that  $\dta(\t)=\om(\t)$ in the tame case \cite{SingKha}.
We give a short proof of this fact and we find an explicit formula for this invariant in terms of the data supported by any MLV chain of $v_F$ (Corollary \ref{lambdas2}).

In Section \ref{secExs}, we display several rank-one examples, including some counterexamples to \cite[Thm. 3.4]{OkS}. Some defect examples suggest that the presence of defect does not increase the difference   $\om(\t)-\dta(\t)$. This leads to the following conjecture.\e

\nn{\bf Conjecture. }\emph{For a rank-one  Henselian $(K,v)$, take $L=K(\t)$ for some separable $\t\in\kb$. If $Lv/K\!v$ is separable and $e(L/K)$ is not divisible by $\chr(K\!v)$, then $\dta(\t)=\om(\t)$}.    

\section{Mac Lane-Vaqui\'e chains and Okutsu frames}\label{secMLVOF}

Recall that $\g:=v\kb =vK\otimes_\Z\Q$ is the divisible hull of $vK$. From now on, we shall write $\gi$ instead of $\g\cup\{\infty\}$.

 \subsection{Valuations on $\kx$}\label{subsecTG}

 A valuation on $\kx$, taking values in $\g$, is a mapping $\,\mu\colon \kx\to \gi$, 
satisfying the following conditions for all $f,g\in \kx$:
\begin{itemize}
 \item \ $\mu(1)=0$, \ $\mu(0)=\infty$, 
 \item \ $\mu(fg)=\mu(f)+\mu(g)$, 
  \item \ $\mu(f+g)\ge\min\{\mu(f),\mu(g)\}$.
  \end{itemize}\e

 Let $\ttt(\g)$ be the set of all valuations on $\kx$, taking values in $\g$, and extending $v$. This set has the structure of  a tree (all intervals are totally ordered), with respect to the partial ordering described in (\ref{ordering}).
 The maximal elements in $\ttt(\g)$ are said to be \textbf{leaves} of this tree.
 
For any $\mu\in\ttt(\g)$, the \textbf{support} of $\mu$ is the prime ideal  
 \[
 \p:=\supp(\mu)=\mu^{-1}(\infty)\in\op{Spec}(\kx). 
 \]
Only the valuations with $\p=0$ can be extended to valuations on the field $K(x)$. Actually, every valuation $\mu$ induces in an obvious way a valuation $\mub$ on the field of fractions of $\kx/\p$.
 We denote by $\gm$ and $\km$ the value group and residue field of $\mu$, which are defined as the value group and residue field of $\mub$, respectively.  
 \e

\nn{\bf Remark. }For instance, if $F\in\irr(K)$ is a generator of $L/K$, then the valuation $v_F$ has support $\p=F\kx$  and the valuation $\overline{v_F}$ on $\kx/(F)\simeq L$ can be identified with $v$. By definition, we have $\g_{v_F}=vL$ and  $k_{v_F}=Lv$.  

\begin{definition}\label{resTrans}
 We say that $\mu$ is \textbf{residue-transcendental} if 
 $\p=0$ and $\km/K\!v$ is transcendental. In this case, the transcendence degree of  $\km/K\!v$ is equal to one.
 \end{definition}

 The {\bf graded algebra} of $\mu$ is the integral domain $\ggm=\bigoplus_{\alpha \in \Gamma_\mu}\mathcal P_\alpha(\mu)/\mathcal P_\alpha^+(\mu)$, where
 \[
 \mathcal P_\alpha(\mu)=\{f\in K[x]\mid \mu(f)\geq \al\}\supseteq\mathcal P_\alpha^+(\mu)=\{f\in K[x]\mid \mu(f)> \al\}.
 \]
 Every $f\in K[x]\setminus\p$ has a homogeneous initial coefficient  $\inm f\in\ggm$, defined as the image of $f$ in $\mathcal P_{\mu(f)}/\mathcal P_{\mu(f)}^+$. 
 
%\begin{lemma}\label{simple}  If $\p\ne0$, then the graded algebra $ \ggm$ is simple; that is, every nonzero homogeneous element is a unit. \end{lemma}

  \begin{definition}\label{keypolyno}
  	Take a monic $\phi\in\kx$ and let $I(\phi)$ be the principal ideal $(\inm\phi)\ggm$. 
  	We say that $\phi$ is \textbf{$\mu$-minimal} if $I(\phi)$ contains no initial coefficient \ $\inm f$ with  $\deg(f)< \deg(\phi)$.	
 	We say that $\phi$ is a \textbf{key polynomial}  for $\mu$ if it is $\mu$-minimal and  $I(\phi)$ is a prime ideal.
 \end{definition}

 \nn{\bf Notation. }Let $\kpm$ be the set of all key polynomials for $\mu$. 
If $\mu<\nu$ in $\ttt(\g)$, then we denote by $\ty(\mu,\nu)$ the set of all monic $\varphi\in\kx$ of minimal degree satisfying $\mu(\varphi)<\nu(\varphi)$. 

We have $\ty(\mu,\nu)\sub\kpm\sub\irr(K)$.

\begin{lemma}\cite[Thm. 1.15]{Vaq}\label{td}
If $\mu<\nu$ in $\ttt(\g)$, then for all $\varphi\in\ty(\mu,\nu)$, $f\in\kx$, we have
 \[
  \mu(f)<\nu(f) \ \sii \ \inm \varphi\ \mbox{ divides }\ \inm f \ \mbox{ in }\ \ggm.
 \]
\end{lemma}

 \begin{proposition}\cite[Thm. 3.9]{KP}\label{weightMu}
 	Let $\phi\in\kpm$. Then,
$ \mu(f)/\deg(f)\le \mu(\phi)/\deg(\phi)$, for all monic $f\in\kx$. 
Equality holds if and only if $f$ is $\mu$-minimal.
 \end{proposition}
 
  The following  criterion  for the existence of key polynomials follows from \cite[Thm 4.4]{KP} and \cite[Thm 2.3]{MLV}.
 
 \begin{theorem}\label{existence}
 	For every $\mu\in\ttt(\g)$ the following conditions are equivalent.
 	\begin{enumerate}
 		\item[(a)] $\kpm\ne\emptyset$.
 		\item[(b)] $\mu$ is residue-transcendental.
 		\item [(c)] $\mu$ is not a leaf of $\ttt(\g)$. 
 	\end{enumerate}
 \end{theorem}

\begin{corollary}\label{leaf}
If $\p\ne0$, then $\kpm=\emptyset$ and $\mu$ is a leaf of $\ttt(\g)$. 
\end{corollary}

\begin{proof}
	We have $\p=G\kx$ for some monic, irreducible $G\in\kx$.
	Since $\kx/\p$ is a field, for every $f\in\kx\setminus\p$  there exists $g\in\kx$ such that $fg\equiv 1\md{\p}$.  Thus, $(\inm f)(\inm g)=\inm 1\,$ in $\ggm$.
	Therefore,  $\kpm=\emptyset$ because $\ggm$  contains no homogeneous prime elements. By Theorem \ref{existence}, $\mu$ is a leaf in $\ttt(\g)$. 
\end{proof}

 \begin{definition}\label{degMu}
  Let $\mu\in\ttt(\g)$ be residue-transcendental. We define the \textbf{degree} of $\mu$ as the minimal degree of a key polynomial for $\mu$. We denote it by $\deg(\mu)$.
\end{definition}

By Theorem \ref{existence}, the residue-transcendental valuations admit augmentations in the tree $\ttt(\g)$. For our purposes, we consider only two specific types of augmentations. 
 
 \begin{definition}\label{defOrd}
For some $\mu\in\ttt(\g)$, take $\phi\in\kpm$ and $\ga\in\gi$ such that $\ga>\mu(\phi)$. The \textbf{ordinary augmentation} $ \nu=[\mu;\,\phi,\ga]$
 is the valuation defined in (\ref{ordinary}). 
 \end{definition}

 If $\ga=\infty$, then $\supp(\nu)=\phi\kx$ and $\nu$ is a leaf of $\ttt(\g)$ by Corollary \ref{leaf}. 
 If $\ga<\infty$, then  $\phi$ becomes a key polynomial  of minimal degree of  $\nu$ \cite[Cor. 7.3]{KP}. Thus, $\deg(\nu)=\deg(\phi)$.\e

 	Let $I$ be an infinite well-ordered set and $\aa=\left(\ri\right)_{i\in I}$ a family of valuations in $\ttt(\g)$ such that $\ri<\rj$ whenever $i<j$ in $I$. 

We say that $g\in\kx$ is \textbf{$\aa$-unstable} if $\ri(g)<\rj(g)$ whenever $i<j$ in $I$.   An $\aa$-unstable polynomial of minimal degree is said to be a \textbf{limit key polynomial} for $\aa$. Let $\kpi(\aa)$ be the set of all limit key polynomials for $\aa$.

 \begin{definition}\label{continuous}
 We say that $\aa$ is a \textbf{continuous family} if all valuations $\ri$ have the same degree and $\kpi(\aa)\ne\emptyset$.
 \end{definition}
   
Take a continuous family  $\aa$ and $\phi\in\kpi(\aa)$. Let $m=\deg(\phi)$. Every $a\in\kx$ with $\deg(a)<m$ is necessarily $\aa$-stable; that is, there exists some $i_0\in I$ such that
\[
\ri(a)=\rj(a)\quad\mbox{ for all }\ j>i\ge i_0.
\]  
Let us denote this stable value by $\rho_\aa(a)$.

 \begin{definition}\label{defLim}
Let $\aa=\left(\ri\right)_{i\in I}$ be a continuous family. Take $\phi\in\kpi(\aa)$ and $\ga\in\gi$ such that $\ga>\ri(\phi)$ for all $i\in I$. The \textbf{limit augmentation} $ \nu=[\aa;\,\phi,\ga]$
	is the valuation in $\ttt(\g)$ defined as follows on $\phi$-expansions:
\[
g=\sum\nolimits_{0\le n}a_n\phi^n,\quad \deg(a_n)<\deg(\phi) \imp
\nu(g)=\min_{0\le n}\{\rho_\aa(a_n)+n\ga\}.
\]
\end{definition}
   
If $\ga=\infty$, then $\supp(\nu)=\phi\kx$ and $\nu$ is a leaf of $\ttt(\g)$ by Corollary \ref{leaf}. 
 If $\ga<\infty$, then $\phi$ becomes a key polynomial  of minimal degree of $\nu$ \cite[Cor. 7.13]{KP}. Thus, $\deg(\nu)=\deg(\phi)$.  

%   such that its minimal valuation is $\mu$
\subsection{MLV chains}\label{subsecMLV}
%Okutsu frames of $F$ were defined in \cite{okutsu} for discrete rank-one valuations, inspired in the work of Okutsu in \cite{Ok}. They were extended to defectless valuations of arbitrary rank in \cite{defless}, and to the general case including defect in \cite{OkF}. 

Let $F\in\kx$ be a generator of $L/K$. A celebrated theorem of Mac Lane-Vaqui\'e \cite{Vaq, MLV} proves the existence of a finite chain of augmentations in $\ttt(\g)$, with end node $v_F$:
\begin{equation}\label{depthMLV}
v\ \stackrel{\phi_0,\ga_0}\lra\  	\mu_0\ \stackrel{\phi_1,\ga_1}\lra\  \mu_1\ \stackrel{\phi_2,\ga_2}\lra\ \cdots
	\ \lra\ \mu_{r-1} 
	\ \stackrel{F,\infty}\lra\ \mu_r=v_F, 
\end{equation}
The formal augmentation $v\to\mu_0$ was described in (\ref{mu0}). It is easy to check that $\phi_0$ is a key polynomial of $\mu_0$; thus, $\deg(\mu_0)=1$. Denote $\phi_r:=F$, $\ga_r:=\infty$.
For $0<i\le r$, each augmentation $\mu_{i-1} 
\to \mu_i$ is of one of the following types:\e

\emph{Ordinary augmentation}: \ $\mu_{i}=[\mu_{i-1};\, \phi_{i},\ga_{i}]$, for some $\phi_{i}\in\kp(\mu_{i-1})$.\e

\emph{Limit augmentation}:  \ $\mu_{i}=[\aa_{i-1};\, \phi_{i},\ga_{i}]$,  for some $\phi_{i}\in\kpi(\aa_{i-1})$, where $\aa_{i-1}$ is a continuous family admitting $\mu_{i-1}$ as its first valuation.\e

By Theorem \ref{existence}, the valuations  $\mu_0,\dots,\mu_{r-1}$ are residue-transcendental. 
Also,
$$
m_i:=\deg(\mu_i)=\deg(\phi_i),\quad \mu_i(\phi_i)=\ga_i,\qquad 0\le i\le r.
$$

\defn
{\it A chain of mixed augmentations as in (\ref{depthMLV}) is said to be  a \textbf{\mlv (MLV) chain} of $v_F$ if it satisfies:}
	\[
	1=m_0<\cdots<m_r=n\qquad\mbox{and}\qquad v_F(\phi_i)=\ga_i \quad\mbox{for all }\,0\le i<r.
	\]
%\begin{itemize}
%	\item $\,1=m_0<\cdots<m_r=n$.
%	\item $\,v_F(\phi_i)=\ga_i$ \ for all $\,0\le i<r$. 
%end{itemize}}\e

Let us recall some intrinsic data of $v_F$ supported by every MLV chain of $v_F$. For every $\mu\in\ttt(\g)$, let $\ka_\mu$ be the relative algebraic closure of $K\!v$ inside $k_\mu$. Also, let $\g_{\mu_{-1}}:=vK$.\e

\nn{\bf Relative ramification indices and  inertia degrees.} For all $0\le i<r$, we define
\[
 e_i:=\left(\g_{\mu_i}\colon \g_{\mu_{i-1}}\right),\qquad f_i:=\left[\ka_{\mu_{i+1}}\colon \ka_{\mu_i}\right].
\]
 
\nn{\bf Relative defect.}
For all $0\le i<r$, we define
\[
d_i:=
\begin{cases}
1,&\mbox{ if }\,\mu_i\to\mu_{i+1}\ \mbox{ is ordinary},\\
\deg(\phi_{i+1})/\deg(\phi_i)	&\mbox{ if }\,\mu_i\to\mu_{i+1}\ \mbox{ is limit}.
\end{cases}
\]

\begin{proposition}\label{efd}\cite{Vaq,MLV}%\mbox{\null}
	
\begin{enumerate}
	\item [(a)] If $\,\mu_i\to\mu_{i+1}$ is a limit augmentation, then $e_i=f_i=1$.
	\item [(b)] For all $i\ge 0$ we have: $\ m_{i+1}=e_if_id_im_i$.
	\item [(c)] $e(L/K)=e_0\cdots e_{r-1}$, \ $f(L/K)=f_0\cdots f_{r-1}$, \ $d(L/K)=d_0\cdots d_{r-1}$.
\end{enumerate}	
\end{proposition}

%We emphasize that this result does not hold if $(K,v)$ is not Henselian.

Now, let us analyze some (non-intrinsic) elements in $\g$ associated to the MLV chain too. 

\begin{proposition}\label{Fminimal}
For every $\,0< i\le r$, the polynomial $\phi_{i}$ is $\mu_{i-1}$-minimal. 
\end{proposition}
	
\begin{proof}
If $\mu_{i}=[\mu_{i-1};\,\phi_{i},\ga_{i}]$ is an ordinary augmentation, then $\phi_{i}$ is a key polynomial for $\mu_{i-1}$; hence, it is $\mu_{i-1}$-minimal. 

Suppose that $\mu_{i}=[\aa_{i-1};\,\phi_{i},\ga_{i}]$ is a limit augmentation. Let $\aa_{i-1}=\left(\ri\right)_{i\in I}$ with $\mu_{i-1}=\rho_{i_0}$, where $i_0=\min(I)$.  Since $\phi_{i}$ is a limit key polynomial for $\aa_{i-1}$, it satisfies:
\[
\mu_{i-1}(\phi_{i})<\ri(\phi_{i})<v_F(\phi_{i})\quad\mbox{for all } \ i\in I,\ i>i_0. 
 \]
Take any $\varphi\in \ty(\mu_{i-1},v_F)$. By Lemma \ref{td}, $\varphi$ is a key polynomial for $\mu_{i-1}$   
and $\inn_{\mu_{i-1}}\varphi$ divides $\inn_{\mu_{i-1}} \phi_{i}$ in the graded algebra $\gg_{\mu_{i-1}}$. In particular, $\inn_{\mu_{i-1}} \phi_{i}$ is not a unit and \cite[Cor. 5.3]{NN} shows that $\phi_i$ is $\mu_{i-1}$-minimal.
\end{proof}\e

For every $i\ge 0$, the Newton polygon $N_{\mu_i,\phi_i}(F)$ is one-sided of slope $\ga_i$ \cite[Thm. 5.2]{NN}. Thus, it is natural to say that $\ga_0,\dots,\ga_{r-1}$ are the \textbf{slopes} of $F$, with respect to this particular MLV chain of $v_F$.

We may consider certain \textbf{secondary slopes} $\la_0,\dots,\la_{r-1}\in\g$ defined as follows:
\[
\la_0:=\ga_0,\qquad \la_i:=\ga_i-\mu_{i-1}(\phi_i)>0\quad\mbox{for all }\ 0<i<r.
\]
The next identity was proven in \cite{defless}, under the assumption that $L/K$ was defectless.

\begin{lemma}\label{lambdas}
	For all $0\le i< r$, we have \  $
\dfrac{\ga_i}{m_i}=\dfrac{\la_0}{m_0}+\cdots+\dfrac{\la_i}{m_i}$.
\end{lemma}

\begin{proof}
This is obvious for $i=0$. Suppose that  $i>0$. By Propositions \ref{weightMu} and \ref{Fminimal}, 
\[
\dfrac{\ga_{i-1}}{m_{i-1}}=\dfrac{\mu_{i-1}(\phi_{i-1})}{m_{i-1}}=\dfrac{\mu_{i-1}(\phi_i)}{m_i}=\dfrac{\ga_i-\la_i}{m_i}.
\] 
This proves the claimed identity by a recursive argument.
\end{proof}\e

The main advantage of MLV chains is that they support data intrinsically associated to the valuation $v_F$, like
the sequence $(m_n)_{n\ge0}$, the character ``ordinary" or ``limit" of each augmentation step $\mu_i\to\mu_{i
+1}$ and the relative data $e_i$, $f_i$, $d_i$, of each step \cite{MLV}.  

In particular, all MLV chains of $v_F$ have the same length $r$. This common length is said to be the \textbf{depth} of $v_F$. We say that $r$ is the depth of $F$ as well. 

These data are not intrinsically associated to the extension $L/K$. They depend on the choice of the generator $F$. In Section \ref{secExs}, we shall see some examples of this dependence.

\subsection{Okutsu frames}\label{subsecOF}
Let us fix some  $F\in\irr(K)$  of degree $n>1$. Recall the weight function (with respect to $F$) that was defined in the Introduction:
\[%\begin{equation}\label{weight}
 w(g):=v_F(g)/\deg(g)\in\g,\quad\mbox{for all }\ g\in\kx\setminus K.
\]%\end{equation}
For $m\in\N$, $m>1$, denote
\[
\kx_m:=\left\{g\in\kx\mid g\mbox{ monic},\ 1\le \deg g<m\right\}.
\]
For every integer $1<m\le n$, consider the set
$$
W_m=W_m(F):=w\left(\kx_m\right)\sub\g.
$$

We say that a  subset $\Phi\sub\irr(K)$ has a ``common degree" if all polynomials in $\Phi$ have the same degree. In this case, we denote this common degree by $\deg\Phi$. 

\begin{definition}\label{defOkF}
An \textbf{Okutsu frame} of $F$ is a finite list \[ \left[\Phi_0,\Phi_1,\dots,\Phi_{r-1},\Phi_r=\{F\}\right]\]
of common-degree subsets of $\irr(K)$,
whose degrees grow strictly:
\[
1=m_0<m_1<\cdots<m_r=n,\qquad m_\ell=\deg\Phi_\ell, \ \ 0\le\ell\le r,
\]
and the following fundamental property is satisfied for all $\,0\le\ell< r$: \e

(OF0) \ For every $g\in\kx_{m_{\ell+1}}$, there exists $\phi\in\Phi_\ell$ such that
$w(g)\le  w(\phi)$.\e

Moreover, the  following additional properties are imposed, for all $0\le\ell< r$: \e

(OF1) \ $\#\Phi_\ell=1$ whenever $\max\left(W_{m_{\ell+1}}\right)$ exists.\e

(OF2) \ If $\max\left(W_{m_{\ell+1}}\right)$ does not exist, we assume that $\Phi_\ell$ is totally ordered with respect to the action of $w$: \ $\phi<\phi'\ \sii \ w(\phi)<w(\phi')$.\e

(OF3) \ For all $\phi\in\Phi_{\ell}$, $\varphi\in\Phi_{\ell+1}$, we have $w(\phi)<w(\varphi)$.
\end{definition}

The following two theorems show that
Okutsu frames of $F$ and MLV chains of $v_F$ are essentially equivalent objects. 

\begin{theorem}\cite[Thm. 4.4]{OkF}\label{MLOk}
Let $F\in\irr(K)$ and consider a MLV chain of $v_F$:
$$
v \ \stackrel{\phi_0,\ga_0}\lra\ \mu_0\ \stackrel{\phi_1,\ga_1}\lra\  \mu_1\ \lra\ \cdots
\ \stackrel{\phi_{r-1},\ga_{r-1}}\lra\ \mu_{r-1}
\  \stackrel{F,\infty}\lra\ \mu_{r}=v_F.
$$
Then, $\left[\Phi_0,\dots,\Phi_{r-1},\Phi_{r}=\{F\}\right]$ is an Okutsu frame of $F$, where for $0\le\ell< r$, we consider the following common-degree subsets  of $\irr(K)$:

{$\bullet$} \ If $\mu_\ell\to\mu_{\ell+1}$ is an ordinary augmentation, we take 	$\Phi_\ell=\left\{\phi_\ell\right\}$.

{$\bullet$} \  If $\mu_\ell\to\mu_{\ell+1}$ is a limit augmentation with  respect to a continuous family $\aa_\ell=\left(\rho_i\right)_{i\in I_\ell}$ of valuations, we take $\Phi_\ell=\left\{\chi_i\mid i\in I_\ell\right\}$ with $\chi_i\in\ty(\rho_i,v_F)$.
\end{theorem}

For any monic polynomial $\phi\in \kx$, the \textbf{truncation} of $v_F$ by $\phi$ is the following function $v_{F,\phi}$, defined on the set $\kx$ in terms of $\phi$-expansions:
\[
 g=\sum_{0\le s}a_s\phi^s,\ \deg(a_s)<\deg(\phi)\ \imp\
v_{F,\phi}(g):=\min_{0\le s}\{v_F(a_s\phi^s)\}. 
\]
This function is not necessarily a valuation on $\kx$ \cite{Dec}.

\begin{theorem}\cite[Thm. 4.5]{OkF}\label{OkML}
Let $\left[\Phi_0,\dots,\Phi_{r-1},\Phi_r=\{F\}\right]$ be an Okutsu frame of $F\in\irr(K)$. For all $0\le \ell\le r$, choose an arbitrary $\,\phi_\ell\in\Phi_\ell$ and denote $\ga_\ell=v_F(\phi_\ell)$. Then, the truncation of $v_F$ by $\phi_\ell$ is a valuation $\mu_\ell$ fitting into a MLV chain of $v_F$:
\[
v \ \stackrel{\phi_0,\ga_0}\lra\ \mu_0\ \stackrel{\phi_1,\ga_1}\lra\  \mu_1\ \lra\ \cdots
\ \stackrel{\phi_{r-1},\ga_{r-1}}\lra\ \mu_{r-1}
\  \stackrel{F,\infty}\lra\ \mu_{r}=v_F.
\]

If $\Phi_\ell=\left\{\phi_\ell\right\}$, then $\mu_{\ell+1}=[\mu_\ell;\,\phi_{\ell+1},\ga_{\ell+1}]$ is an ordinary augmentation.

If $\Phi_\ell=\left\{\chi_i\mid i\in J_\ell\right\}$, then $\mu_{\ell+1}=[\aa_\ell;\,\phi_{\ell+1},\ga_{\ell+1}]$ is a limit augmentation with res\-pect to the essential continuous family $\aa_\ell=\{\mu_\ell\}\cup\left(\rho_i\right)_{i\in I_\ell}$, where $\rho_i=[\mu_\ell;\,\chi_i,v_F(\chi_i)]$ and $I_\ell\sub J_\ell$ contains the indices  $i$ such that $v_F(\chi_i)>\ga_\ell$.
\end{theorem}

\begin{corollary}\label{CorOF}
Let $\left[\Phi_0,\dots,\Phi_{r-1},\Phi_r=\{F\}\right]$ be an Okutsu frame of $F\in\irr(K)$. Then:
\begin{enumerate}
\item $\dep(F)=r$.
\item $1=m_0\mid m_1\mid\cdots\mid m_{r-1}\mid \deg(F)$.
\item $F$ is defectless if and only if $\,\#\Phi_\ell=1$ for all $0\le\ell\le r$.  
\item The set $\,\Phi_0\cup\cdots\cup\Phi_{r-1}\cup\{F\}$ is a complete set of abstract key polynomials for $v_F$.
\end{enumerate}
\end{corollary}

\begin{proof}
Items (1), (2) and (3) follow from the results in Section \ref{subsecMLV}.

Item (4) was proved in \cite[Thm. 4.8]{OkF}.
\end{proof}\e

For the definition of complete sets of abstract key polynomials, see \cite{Dec} and \cite{NS2018}.

\section{Weight and distance of polynomials}\label{secWD}

For any $\al\in\kb$ we denote by $\irr_K(\al)\in \irr(K)$ its minimal polynomial over $K$. 
We fix, once and for all, some $\t\in\kb$ and denote 
$n:=\dgk \t>1$, $F:=\irr_K(\t)$.

\begin{definition}
	For all monic $f\in\kx$, we define its  \textbf{distance} of $f$ to $\t$ as:
	\[
	d(f):=\max\left\{v(\t-\al)\mid \al\in Z(f)\right\}\in\gi.
	\]
\end{definition}

Clearly, the weight 
\[
w(f)=v_F(f)/\deg(f)=v(f(\t))/\deg(f)
\]
is simply the average of the values $v(\t-\al)$ for $\al$ running in the multiset $Z(f)$. Thus, 
\[
w(f)\le d(f)\quad\mbox{for all}\quad f\in\kx.
\]

The following  observation is trivial.

\begin{lemma}\label{product}
	If $f=f_1\cdots f_t$ is a product of monic polynomials in $\kx$, then:
	\[
	w(f)\le \max\{w(f_i)\mid 1\le i\le t\},\qquad 
	d(f)= \max\{d(f_i)\mid 1\le i\le t\}.
	\]
\end{lemma}

For all $f\in\irr(K)$, let $\ep(f):=v(\be)$ be the common value of all $\be\in Z(f)$. Clearly, $\ep(f)\le\om(f)$, where $\om(f)$ is Krasner's constant, defined in the Introduction.

\begin{lemma}\label{inseparable}
	Let $f\in\irr(K)$ be inseparable of degree $m$ and take any $\rho\in \g$. 
	By taking $\pi\in K^*$ with $v(\pi)$ sufficiently large, the separable polynomial $\fsep:=f+\pi x$ has degree $m$ and, for a proper ordering of the roots of $f$ and $\fsep$
	\[
	Z(f)=\{\be_1,\dots,\be_m\},\quad Z(\fsep)=\{\al_1,\dots,\al_m\},
	\]
	we have $v(\be_i-\al_i)>\rho$, for all $1\le i\le m$. 
\end{lemma}

\begin{proof}
	Since $f$ is inseparable, we have $m>1$ and $\deg \fsep =m$. Since $f$ is irreducible, we have $f'=0$, so that $\fsep'=\pi\ne0$ and $\fsep$ is separable.  
	
	Denote $\ep:=\ep(f)$, $\om:=\om(f)$, for simplicity.  We may assume that $\rho\ge\om\ge \ep$. 
	
	Take any
	$\pi\in K$ with   $v(\pi)>m\rho-\ep$.
	The Newton polygon $N_{v,x}(f)$ has only one side of slope $-\ep$ and end points $(0,m\ep)$, $(m,0)$.  Since $v(\pi)>(m-1)\ep$, the polynomial $\fsep$ has the same Newton polygon: $N_{v,x}(f)=N_{v,x}(\fsep)$.

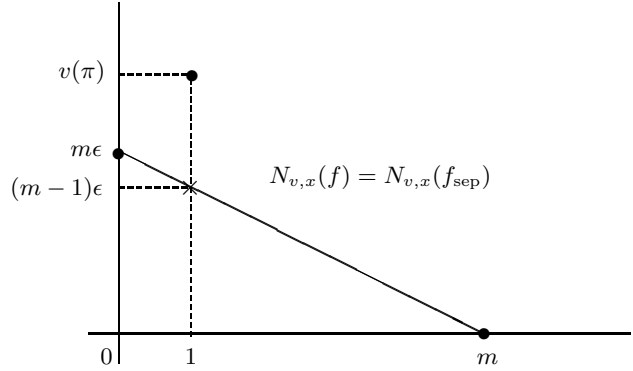
	\begin{center}
		\setlength{\unitlength}{4mm}
		\begin{picture}(19,12.5)
			\put(2.2,9.35){$\bullet$}
			\put(-0.25,6.75){$\bullet$}\put(11.9,0.75){$\bullet$}
			\put(2,5.6){$\times$}
			\put(-1,1){\line(1,0){18}}
			\put(0,0){\line(0,1){12}}
			\put(0,7.07){\line(2,-1){12}}\put(0,7.1){\line(2,-1){12}}
			\multiput(2.4,.9)(0,.3){29}{\vrule height2pt}
			\put(5,6){\begin{footnotesize}$N_{v,x}(f)=N_{v,x}(\fsep)$\end{footnotesize}}
			\put(-1.6,6.9){\begin{footnotesize}$m\ep$\end{footnotesize}}
			\put(-3.6,5.6){\begin{footnotesize}$(m-1)\ep$\end{footnotesize}}
			\put(-2,9.5){\begin{footnotesize}$v(\pi)$\end{footnotesize}}
			\put(11.9,0){\begin{footnotesize}$m$\end{footnotesize}}
			\put(2.2,0){\begin{footnotesize}$1$\end{footnotesize}}
			\put(-.6,0){\begin{footnotesize}$0$\end{footnotesize}}
			\multiput(0,9.6)(.3,0){8}{\hbox to 2pt{\hrulefill }}
			\multiput(0,5.85)(.3,0){8}{\hbox to 2pt{\hrulefill }}
		\end{picture}
	\end{center}\e
	
Thus, $v(\al)=\ep$ for all $\al\in Z(\fsep)$ as well.
	Now, for any $\al\in Z(\fsep)$, the inequality
	\[
	\sum\nolimits_{\be\in Z(f)}v(\be-\al)=v(f(\al))=v(\pi\al)>m\rho,
	\]
	implies the existence of some $\be\in Z(f)$ such that $v(\be-\al)>\rho$. 
	Since $\rho\ge \om$,
	for every $\be'\in Z(f)$, $\be'\ne\be$, we have $v(\be'-\al)=v(\be'-\be)\le \om$. Thus, $\al$ is ``very close" to a unique root $\be$ of $f$ (disregarding multiplicities). Applying this argument to all roots of $\fsep$, we see that 
	$v(\be_i-\al_i)>\rho$, for all $1\le i\le m$, for an opportune ordering  of $Z(f)$ and $Z(\fsep)$. 
\end{proof}\e

Let us emphasize that $\fsep$ need not be irreducible. For instance, suppose that $f$ is purely inseparable and let $L$ be the separable extension of $K$ obtained by adjoining all roots of $\fsep$. Then, $f$ is still irreducible over $L[x]$, while $\fsep$ splits completely.

The following result is an immediate consequence of Lemma  \ref{inseparable}.

\begin{corollary}\label{approx}
	Let $\be\in\kb$ and take $\rho\in\g$. There exists a separable $\al\in\kb$ such that
	\[
	\deg_K\al\le\deg_K\be\quad\mbox{ and }\quad v(\be-\al)>\rho. 
	\]	
\end{corollary}

\begin{corollary}\label{samewd}
	Let $f\in\kx$ be a monic polynomial such that $f(\t)\ne0$. Then, there exists a separable monic polynomial $\fsep\in\kx$ such that
	\[
	\deg f=\deg \fsep,\qquad w(f)=w(\fsep)\quad\mbox{ and }\quad 	d(f)=d(\fsep).
	\]
\end{corollary}

\begin{proof}
	If $f$ is separable, we may take $\fsep=f$.
	
	Suppose that $f$ is  inseparable and irreducible of degree $m$. Take $\rho\in\g$ such that $\rho\ge\max\{\om(f),d(f)\}$. Let $\fsep=f+\pi x$, with $\pi\in K^*$ satisfying
	\[
	v(\pi)>m\rho-\ep(f)\quad\mbox{ and }\quad v(\pi)>v(f(\t))-v(\t).
	\]
	
	Since $v(\pi\t)>v(f(\t))$, we have $v(f(\t))=v(\fsep(\t))$, so that $w(f)=w(\fsep)$. As shown in the proof of Lemma \ref{inseparable},
	for an opportune ordering  of the multiset $Z(f)=\{\be_1,\dots,\be_m\}$ and the set $Z(\fsep)=\{\al_1,\dots,\al_m\}$, we have
	\[
	v(\be_i-\al_i)>\rho\ge d(f)\ge v(\t-\be_i), \quad\mbox{for all }1\le i\le m.
	\]
	This implies $v(\t-\be_i)=v(\t-\al_i)$ for all $i$, so that $d(f)=d(\fsep)$.
	
	Finally, suppose that $f=\phi_1\cdots \phi_t$ is a product of monic irreducible polynomials in $\kx$. For each irreducible factor $\phi$ of $f$, take a separable $\phsep$ of the same degree 
	such that $w(\phi)=w(\phsep)$ and $d(\phi)=d(\phsep)$, as indicated above. Let $\fsep$ be the product of all $\phsep$. Clearly, $\deg f=\deg \fsep$.
	
	Since $\deg \phi=\deg \phsep$, we have $v(\phi(\t))=v(\phsep(\t))$ for all $\phi$; hence, $v(f(\t))=v(\fsep(\t))$ and $w(f)=w(\fsep)$.
	Also, 
	\[
	d(f)=\max_{\phi\mid f}\{d(\phi)\}=\max_{\phi\mid f}\{d(\phsep)\}=d(\fsep).
	\]
	This ends the proof.
\end{proof}\e

\nn{\bf Notation. }For $m\in\N$, $m>1$, denote
\[
\kx^{\sep}_m:=\left\{g\in\kx\mid g\mbox{ monic, separable},\ 1\le \deg g<m\right\}.
\]
\[
\irr(K)^{\sep}_m:=\left\{g\in\irr(K)\mid g\mbox{ separable},\ 1\le \deg g<m\right\}.
\]

\begin{corollary}\label{SepCofinal}
	Let $m\in\N$, $1<m\le n$. Then, $d\left(\irr(K)^{\sep}_m\right)=d\left(\kx_m\right)$ and  $w\left(\irr(K)^{\sep}_m\right)$ is cofinal in $w\left(\kx_m\right)$.
\end{corollary}

\begin{proof}
	By Corollary \ref{samewd}, we have 
	\[
	d\left(\kx^{\sep}_m\right)=d\left(\kx_m\right),\qquad 
	w\left(\kx^{\sep}_m\right)=w\left(\kx_m\right).
	\]	
	By Lemma \ref{product}, $d\left(\irr(K)^{\sep}_m\right)=d\left(\kx^{\sep}_m\right)$ and  $w\left(\irr(K)^{\sep}_m\right)$ is cofinal in $w\left(\kx^{\sep}_m\right)$. 
\end{proof}\e

The following result is inspired in \cite[Thm. 2.3]{OkS}.

\begin{proposition}\label{ftalEquiv}
	Let $f,g\in\kx$ be monic, irreducible and separable polynomials. Then,
	$d(f)\le d(g)$ if and only if $w(f)\le w(g)$.
\end{proposition}

\begin{proof}
	Clearly, for any $f\in\kx$ we have
	\[
	d(f)=\infty\ \sii\ f(\t)=0\ \sii\ w(f)=\infty.
	\]
	Thus, the case in which $f(\t)g(\t)=0$ is trivial. Let us assume $f(\t)g(\t)\ne0$.
	
	Let $\rho=\max\{d(f),d(g)\}$. By Corollary \ref{approx}, there exists a separable $\t_{\op{sep}}\in\kb$ such that $v(\t-\t_{\op{sep}})>\rho$. This implies
	\[
	v(\t-\al)=v(\t_{\op{sep}}-\al),\quad v(\t-\be)=v(\t_{\op{sep}}-\be),\quad \forall\al\in Z(f),\ \forall\be\in Z(g).
	\]
	Hence, we may assume that $\t$ is separable. Let $L\supset K$ be a finite Galois extension containing $\t$ and all roots of $f$ and $g$. Let $G=\op{Gal}(L/K)$.
	Choose $\al\in Z(f)$, $\be\in Z(g)$ such that $v(\t-\al)=d(f)$, $v(\t-\be)=d(g)$. For all $\sg\in G$ we have
	\begin{equation}\label{rmk}
		v(\t-\sg(\t))\ge \min\{v(\t-\sg(\al)),v(\sg(\al)-\sg(\t))\}=	
		v(\t-\sg(\al)),
	\end{equation}
	and similarly, $v(\t-\sg(\t))\ge v(\t-\sg(\be))$.
	
	Suppose that $d(f)\le d(g)$; that is, $v(\t-\al)\le v(\t-\be)$. We claim that
	\begin{equation}\label{eachOne}
		v(\t-\sg(\al))\le v(\t-\sg(\be))\quad\mbox{for all}\quad \sg\in G.
	\end{equation}
	Indeed, if $v(\t-\sg(\be))=d(g)$, this follows from $v(\t-\sg(\al))\le d(f)\le d(g)$. On the other hand, if  $v(\t-\sg(\be))<d(g)$, then (\ref{eachOne}) follows from (\ref{rmk}). Indeed,
	\[
	v(\t-\sg(\al))\le v(\t-\sg(\t))=\min\{v(\t-\sg(\be)),v(\sg(\be)-\sg(\t))\}=v(\t-\sg(\be)),
	\]
	because $v(\t-\sg(\be))<d(g)=v(\be-\t)=v(\sg(\be)-\sg(\t))$. 
	
	Since $f$ and $g$  are irreducible, $G$ acts transitively on each of their sets of roots. Hence, we may deduce from (\ref{eachOne}), the desired inequality $w(f)\le w(g)$ as follows:
	\begin{equation}\label{ww}
		\dfrac{\#G}{\deg f}\,v(f(\t))=\sum_{\sg\in G}v(\t-\sg(\al))\le
		\sum_{\sg\in G}v(\t-\sg(\be))=\dfrac{\#G}{\deg g}\,v(g(\t)). 
	\end{equation}
	
	Conversely, suppose $w(f)\le w(g)$. We want to deduce that $v(\t-\al)\le v(\t-\be)$. Suppose that $v(\t-\al)> v(\t-\be)$. Then, the above argument shows that (\ref{ww}) would  hold, but exchanging the role of $f$ and $g$. Moreover, the inequality
	\[
	\sum_{\sg\in G}v(\t-\sg(\be))\le
	\sum_{\sg\in G}v(\t-\sg(\al))
	\]
	would be strict because it is strict for $\sg=1$. This would imply $w(g)<w(f)$, against our assumption.
\end{proof}

\section{Okutsu sequences}\label{secOS}

We keep dealing with some fixed $\t\in\kb$ of degree $n>1$ over $K$ and minimal polynomial $F:=\irr_K(\t)$. 
 For every integer $1<m\le n$, consider the set
\[
V_m=V_m(\t):=\left\{v(\t-\be)\mid \be\in\kb,\ \deg_K\be<m\right\}\sub\g.
\]

We say that a  subset $A\sub\kb$ has a ``common degree" if all its elements have the same degree over $K$. In this case, we shall denote this common degree by $\deg_K A$. \e

\begin{definition}\label{defOkS}
	An \textbf{Okutsu sequence} of $\t$ is a finite list  \[\left[A_0,A_1,\dots,A_{r-1},A_r=\{\t\}\right]\] 
of common-degree subsets of $\kb$ whose degrees grow strictly:
\[
		1=m_0<m_1<\cdots<m_r=n,\qquad m_\ell=\deg_K A_\ell, \ \ 0\le\ell\le r,
\]
	and the following fundamental property is satisfied for all $\,0\le\ell< r$: \e
	
	(OS0) \ For all $\be\in\kb$ with $\deg_K\be<m_{\ell+1}$, we have $v(\t-\be)\le  v(\t-\al)$ for some $\al\in A_\ell$.\e
	
	Moreover, the  following additional properties are imposed, for all $0\le\ell< r$: \e
	
	(OS1) \ $\#A_\ell=1$ whenever $\max\left(V_{m_{\ell+1}}\right)$ exists.\e
	
	(OS2) \ If $\max\left(V_{m_{\ell+1}}\right)$ does not exist, we assume that $A_\ell$ is totally ordered with respect to the following ordering: \ $\al<\al'\ \sii \ v(\t-\al)<v(\t-\al')$.\e
	
	(OS3) \ For all $\al\in A_{\ell}$, $\ga\in A_{\ell+1}$, we have $v(\t-\al)<v(\t-\ga)$.
\end{definition}

\begin{lemma}\label{recurrence}
Let $\left[A_0,A_1,\dots,A_{r-1},A_r=\{\t\}\right]$ be an Okutsu sequence of $\t$ of length $r>1$. For any $\al\in A_{r-1}$ the list  $\left[A_0,A_1,\dots,A_{r-2},\{\al\}\right]$ is an Okutsu sequence of $\al$.
\end{lemma}

\begin{proof}
Since $r>1$, we have $\al\not\in K$. In order to prove the lemma, it suffices to show that for every $\be\in\kb$ such that $\deg_K\be<\deg_K \al$, we have $v(\al-\be)=v(\t-\be)$.
This follows from the conditions (OS0) and (OS3). Indeed, there exists $\al'\in A_{r-2}$ such that  
\[
 v(\t-\be)\le v(\t-\al')< v(\t-\al).
\]
Hence, $v(\al-\be)=v(\t-\be)$.
\end{proof}\e

It is easy to check the existence of Okutsu sequences for all $\t\in \kb$.  Take $A$ of minimal common degree $m$ with the property that the family $\left(v(\t-\al)\right)_{\al\in A}$ is cofinal in $V_n$. Then, consider $A=A_{r-1}$ and $m=m_{r-1}$. We may take then $\al\in A$ such that $v(\t-a)>v(\t-\be)$ for all $\be\in\kb$ with $\deg_K\be<\deg_K\al$ and iterate the argument with $\al$.

The following observation is an immediate consequence of Corollary \ref{approx}.

\begin{lemma}\label{allSep}
Every $\t\in\kb$ admits Okutsu sequences $\left[A_0,A_1,\dots,A_{r-1},A_r=\{\t\}\right]$ such that all elements in $A_0\cup\cdots\cup A_{r-1}$ are separable over $K$.
\end{lemma}

%\subsection{Comparison with MLV chains}

Our aim in this section is to prove the following theorem. 

\begin{theorem}\label{SF}
Let $\t\in\kb$ and $F=\irr_K(\t)$.	Let $\left[A_0,A_1,\dots,A_{r-1},A_r=\{\t\}\right]$ be an Okutsu sequence of $\t$ and consider the subsets 
\[
\Phi_\ell=\left\{\irr_K(\al)\mid \al\in A_\ell\right\}, \quad 0\le\ell\le r.
\]
Then, $\left[\Phi_0,\dots,\Phi_{r-1},\Phi_{r}=\{F\}\right]$ is an Okutsu frame of $F$.

Conversely, let $\left[\Phi_0,\dots,\Phi_{r-1},\Phi_{r}=\{F\}\right]$ be an Okutsu frame of $F$.	For each $0\le \ell\le r$, if $\Phi_\ell=\left(\phi_i\right)_{i\in I_\ell}$, then  take $A_\ell=\left(\al_i\right)_{i\in I_\ell}$, where  $\,\al_i\in Z(\phi_i)$ is chosen so that
\[
v(\t-\al_i)=d(\phi_i)\quad\mbox{for all}\quad i\in I_\ell.
\]	
Then, $\left[A_0,A_1,\dots,A_{r-1},A_r=\{\t\}\right]$ is an Okutsu sequence of $\t$.
\end{theorem}

%\begin{corollary}\label{depOkS}	The length of any Okutsu sequence of $\t$ is equal to $\dep(F)$.\end{corollary}

\nn{\bf Notation. }Given a subset $S\sub\g$ and $\rho\in\g$, we write $S<\rho$ to indicate that $\rho$ is greater than all elements in $S$.  \e

The proof of Theorem \ref{SF} relies in the next two lemmas.

\begin{lemma}\label{maxExists}
	Let $m\in\N$, $\phi\in\irr(K)$ such that $1\le \deg\phi<m\le n$.
	Then, the following conditions are equivalent:
	
	\,\,(i) \ $d(\phi)=\max(V_m)$ and $\phi$ has minimal degree with this property.
	
	(ii) \ $w(\phi)=\max(W_m)$ and $\phi$ has minimal degree with this property.
\end{lemma}

\begin{proof}
	We claim that it suffices to prove the lemma under the assumption that $\phi$ is irreducible and separable over $K$. Indeed, by Corollary \ref{samewd}, there exists a separable monic $\phsep\in\kx$ such that $\deg\phsep=\deg \phi$, $d(\phsep)=d(\phi)$ and $w(\phsep)=w(\phi)$.

By Lemma \ref{product}, we have $d(\phsep)=d(\eta)$, $w(\phsep)\le w(\xi)$, for some irreducible factors $\eta,\xi$ of $\phsep$ in $\kx$. Now, $\deg\eta<\deg\phsep$ contradicts (i), while  	$\deg\xi<\deg\phsep$ contradicts (ii).
Thus, if either (i) or (ii) hold, then $\phsep$ is irreducible. 
	This ends the proof of our claim.

Let us assume that $\phi$ is irreducible and separable.
Since $d(\kx_m)$ is cofinal in $V_m$, the conditions (i) and (ii) can be rewritten as:\e
	
	\,(i) \ \,$d(g)\le d(\phi)$ for all $g\in\kx_m$,\qquad 
\ \,$d(g)< d(\phi)$ for all $g\in\kx_{\deg\phi}$.

	(ii) \ $w(g)\le w(\phi)$ for all $g\in\kx_m$,\qquad 
$w(g)< w(\phi)$ for all $g\in\kx_{\deg\phi}$.\e

By Corollary \ref{SepCofinal}, we can assume that $g$ is irreducible and separable too. Then, the equivalence between (i) and (ii) follows from Proposition \ref{ftalEquiv}.
\end{proof}\e

\begin{lemma}\label{maxNoExists}
Let $m\in\N$ and $\Phi\sub \irr(K)$ be a common-degree family such that $1\le \deg\Phi<m\le n$. Then, the following conditions are equivalent:

\,(i) \ The set $d\left(\Phi\right)$ is cofinal in $V_m$ and  $d\left(\kx_{\deg\Phi}\right)<d(\phi)$ for some $\phi\in \Phi$.

(ii) \ The set $w\left(\Phi\right)$ is cofinal in $W_m$ and $w\left(\kx_{\deg\Phi}\right)<w(\phi)$ for some $\phi\in \Phi$.
 \end{lemma}

\begin{proof}
Arguing as in the proof of Lemma \ref{maxExists}, we can assume that all $\phi\in\Phi$ are separable. Then, by completely analogous arguments, the result follows from Corollary \ref{SepCofinal} and Proposition \ref{ftalEquiv}.
 \end{proof}\e

\nn{\bf Proof of Theorem \ref{SF}. }Take $B_0:=\{\t\}$, $n_0:=n$, and $\La_0:=\{F\}$.
If $V_n$ admits a maximal element, then we take $\al\in\kb$ with $\deg_K\al$ minimal such that $v(\t-\al)=\max(V_n)$. Consider 
\[
B_1=\{\al\},\quad n_1=\deg_K\al;\qquad \La_1=\{\phi\},\ \phi=\irr_K(\al).
\]
By Lemma \ref{maxExists}, $w(\phi)=\max(W_n)$ and $\phi$ has minimal degree among all monic polynomials with this property.

If $V_n$ has no maximal element, we take a common-degree family $B_1=\left(\al_i\right)_{\ii}$ such  that the values  $v(\t-\al_i)$ are cofinal in $V_n$ and $\deg_K B_1$ is minimal among all common-degree families with this property. Consider  
\[
n_1=\deg_K B_1;\qquad \La_1=\left(\phi_i\right)_{\ii},\ \phi_i=\irr_K(\al_i) \mbox{ for all }\ii.
\]
By Lemma \ref{maxNoExists}, $w(\La_1)$ is cofinal in $W_n$ and $\deg\La_1=n_1$ is minimal among all common-degree families with this property.

A similar procedure with the sets $V_{n_1}$ and $W_{n_1}$ leads to similar common-degree sets $B_2$, $\La_2$ of degree $n_2<n_1$. 
The iteration of this procedure ends with some $B_r$, $\La_r$ of degree one. Then, obviously the list $[A_0,\dots,A_{r-1},A_r=\{\t\}]$ obtained by reversing the ordering: 
\[
A_i:=B_{r-i},\qquad m_i:=n_{r-i}\ \mbox{ for all }\ 0\le i\le r,
\] is an Okutsu sequence of $\t$. By Lemmas \ref{maxExists} and \ref{maxNoExists}, the list $[\Phi_0,\dots,\Phi_{r-},\Phi_r=\{F\}]$ where $\Phi_i:=\La_{r-i}$ for all $i$, is an Okutsu frame of $F$. This ends the proof of Theorem \ref{SF}.\hfill{$\Box$}\e

The following remark is a direct consequence of Theorems \ref{MLOk}, \ref{OkML} and \ref{SF}.

\begin{corollary}
 The length  of every Okutsu sequence $[A_0\dots,A_r]$ of $\t$ is equal to the depth of $F$. Moreover, the set $A_i$ is a one-element set if and only if the augmentation $\mu_i\to\mu_{i+1}$ is ordinary in every MLV chain of $v_F$.
\end{corollary}

%\subsection{Main invariant and Krasner's constant}

\section{Main invariant and Krasner's constant}\label{secTame}

\begin{definition}
 A separable and defectless extension $L/K$ is said to be \textbf{tame} if $Lv/K\!v$  is separable and $e(L/K)$ is not divisible by $\chr(K\!v)$.
\end{definition}

Let $K^s$ be the separable closure of $K$ in $\kb$.
The \textbf{ramification subgroup} of $\gal(K^s/K)$ is defined as follows:
\[
G\ram=\left\{\sigma\in\gal(K^s/K)\mid v(\sigma(c)-c)>v(c),\,\forall\, c\in (K^s)^*\right\}.
\]
Its fixed field $K\ram=(K^s)^{G\ram}$ is called the \textbf{ramification field} of the extension $K^s/K$. This field is the unique maximal tame extension of $K$ in $\kb$. More precisely, for every algebraic extension $L/K$, the subfield 
\[
L\tame:= L\cap K\ram
\]
is the unique maximal tame subextension of $L/K$.\e

Let $\t\in K^s$ and $F=\irr_K(\t)$. We saw in the Introduction that 
\[
\dta(\t) \le \om(\t),
\]
where $\dta(\t)$ is the main invariant of $\t$ and $\om(\t)$ is its Krasner's constant.  If $K(\t)/K$ is tame, then $\dta(\t)=\om(\t)$ \cite{SingKha}, and there are several formulas computing $\dta(\t)$ in  terms of different invariants \cite{BrownMerzel,BrownMerzel2,PPZ}.

In this Section, we focus our attention on an apparently stronger result from \cite{OkS} which, unfortunately, is false.\e

\nn{\bf False statement. }\cite[Thm. 3.4]{OkS} \emph{Let $[\al_0,\dots,\al_{r-1},\t]$ be an Okutsu sequence of a separable and defectless $\t\in\kb$. If $K(\al_{r-1})/K$ is tame, then $\dta(\t)=\om(\t)$.}\e 

We shall reproduce the ideas of \cite{OkS}, obtaining a short proof of the equality $\dta(\t)=\om(\t)$ in the tame case. Also, we shall save from \cite{OkS} a nice computation of $\dta(\t)$ in terms of the secondary slopes of any MLV chain of $v_F$.\e

Throught this Section, $\t\in\kb$ will be assumed to be separable and defectless over $K$.  By Theorem \ref{SF} and Corollary \ref{CorOF}, we may fix a classical Okutsu sequence
\[
 [\al_0,\al_1,\dots,\al_{r-1},\al_r=\t],
\]
of degrees $1=m_0\mid m_1\mid\cdots\mid m_{r-1}\mid m_r=n$.
We denote
$$
\dta_{-1}:=-\infty<\delta_0:=v(\t-\al_0)<\dots<\delta_{r-1}:=v(\t-\al_{r-1})<\delta_r:=\infty.
$$
Note that $\dta_{r-1}=\dta(\t)$.
By Lemma \ref{allSep}, we may assume that $\al_0,\dots,\al_{r-1}$ are separable over $K$. By Lemma \ref{recurrence}, all $\al_i$ are defectless.

%Note that  $\delta_i=v(\al_{i+1}-\al_i)=\delta_K(\al_{i+1})$, for all $0\le i\le r$.

The next (crucial) result is inspired in the original ideas of Okutsu \cite{okutsu,Ok}.

\begin{proposition}\label{Hi} For some index $0\le i\le r$, let  $M/K$ be a finite Galois extension containing $K(\t,\al_i)$. Let $G=\gal(M/K)$ and consider the subgroups
$$
H_i=\{ \sigma \in G \mid v(\t- \sigma(\t))> \dta_{i-1})\} \supseteq \overline{H}_i=\{ \sigma \in G \mid v(\t- \sigma(\t))\geq \dta_i\}.
$$
Let $M^{H_i} \subset M^{\overline{H}_i}\subset M$ be the respective fixed fields. Then,
$$\mbox{\rm $K(\al_i)\tame$}\sub M^{H_i}\sub M^{\overline{H}_i}\sub K(\t)\cap K(\al_i).$$
\end{proposition}

\begin{proof}
In order to see that $M^{\overline{H}_i}\subset K(\t)\cap K(\al_i)$, it suffices to show that all  $\sigma \in G$ fixing $\t$ or $\al_i$ belong to $\overline{H}_i$. If $\sigma(\t)=\t$, then obviously $\sigma \in \overline{H}_i$. If $\sigma(\al_i)=\al_i$, then 
\[
v(\sigma(\t)-\al_i)=v(\sigma(\t)-\sigma(\al_i))=v(\t-\al_i)= \dta_i. 
\]
Hence, $\sg\in \overline{H}_i$, because
$v(\t-\sigma(\t))\geq \min \{v(\t-\al_i),\,v(\al_i-\sigma(\t))\}=\dta_i$.

Finally let us prove that $K(\al_i)\tame\subset M^{H_i}$. We must prove that 
\begin{equation}\label{aim}
H_i \sub \{\sigma \in G \mid v(\sigma(c)-c)>v(c),\quad \forall c \in K(\al_i)^*\}.
\end{equation}
If $i=0$, then $\al_0\in K$ and (\ref{aim}) is obvious. If $i>0$, then every $c \in K(\al_i)^*$ can be written as $c=g(\al_i)$ for some $g \in K[x]$ with  $\deg(g)<m_i$. 
Write $g= a\prod_{\xi\in\z(g)} (x-\xi)$. For all $\sg\in G$ we have
\begin{equation}\label{cc}
\dfrac{\sg(c)}c=\dfrac{g(\sigma(\al_i))}{g(\al_i)}= \prod_\xi \dfrac{\sigma(\al_i)-\xi}{\al_i-\xi} = \prod_\xi \left(1+ \dfrac{\sigma(\al_i)-\al_i}{\al_i-\xi} \right).	
\end{equation}

If $\sg \in H_i$, then $v\left(\sigma(\t)-\t\right) > \dta_{i-1}$, so that
\[
v\left(\sigma(\al_i)-\al_i\right)\geq \min \left\{v\left(\sigma(\al_i)-\sigma(\t)\right),\,v\left(\sigma(\t)-\t\right),\, v\left(\t-\al_i\right)\right\} > \dta_{i-1}
.\] 
By the minimality of $m_i$,  we have $v(\t-\xi)\leq \dta_{i-1}$ for every $\xi\in\z(g)$. Hence,
\[
v(\al_i-\xi)=\min\left\{v(\al_i-\t),\,v(\t-\xi)\right\}=v(\t-\xi)\le \delta_{i-1},
\]
and this implies $v((\sigma(\al_i)-\al_i)/(\al_i-\xi)) >0$. By (\ref{cc}), we deduce 
$v \left(\dfrac{\sigma(c)}{c}-1\right)>0$, and this proves (\ref{aim}).
\end{proof}\e

\begin{corollary}\label{equality}
If for some $0\le i\le r$  the element $\al_i$ is tame, then $H_i=\overline{H}_i$. In particular, if $\t$  is tame, then $\dta(\t)=\om(\t)$.
\end{corollary}

\begin{proof}
If $\al_i$ is tame, then all inclusions in the chain  
	\[
K(\al_i)=	\mbox{\rm $K(\al_i)\tame$}\sub M^{H_i}\sub M^{\overline{H}_i}\sub K(\al_i)
	\]
must be equalities. By Galois theory, we have $H_i=\overline{H}_i$. 
For $i=r$, we have $\al_r=\t$ and $\dta_r=\infty$. The equality $H_r=\overline{H}_r$ implies  $\dta(\t)=\om(\t)$.
\end{proof}

\begin{corollary}\label{alltame}
If $\t$ is tame, then $\al_0,\dots,\al_{r-1}$ are tame and
\begin{equation}\label{chainAlpha}	
K=K(\al_0)\sub K(\al_1)\sub\cdots\sub K(\al_{r-1})\sub K(\t). 
\end{equation}
\end{corollary}

\begin{proof}
	Denote $K_i=K(\al_i)$ for all $0< i<r$.
By \cite[Lem. 2.9]{OkS}, the Okutsu sequence of $\t$ is a distinguished chain. In the terminology of \cite{AK1}, the ``fundamental principle" states that we have two chains of inclusions:
\[
vK\sub vK_1\sub\cdots\sub vK_{r-1}\sub vL,\qquad 
K\!v\sub (K_1)v\sub\cdots\sub (K_{r-1})v\sub Lv.
\]
If $\t$ is tame, then all extensions $(K_{i})v/K\!v$ are separable and all ramification indices $e(K_{i}/K)$ are not divisible by $\chr(K\!v)$. Therefore, $\al_0,\dots,\al_{r-1}$ are tame. The chain of inclusions (\ref{chainAlpha}) follows from 
Proposition \ref{Hi}.
\end{proof}\e

In the next result, we denote by $\left\{\delta_0^{t_0},\dots,\delta_{r-1}^{t_{r-1}}\right\}$ the multiset whose underlying set is $\left\{\delta_0,\dots,\delta_{r-1}\right\}$ and each element $\dta_i$ appears with multiplicity $t_i$.

\begin{theorem}\label{Valuesti}
If $\t\in\kb$  is tame, then the following multisets of cardinality $n-1$ coincide:
$$
\left\{v\left(\t-\t'\right) \mid \t'\in \op{Z}(F), \ \t'\ne\t\right\}=\left\{\delta_0^{t_0},\dots,\delta_{r-1}^{t_{r-1}}\right\},
$$
where \ $t_i=(n/m_i)-(n/m_{i+1})$ \ for all $0\le i< r$.
\end{theorem}

\begin{proof} Let $M/K$ be a finite Galois extension of $K$ containing $K(\t)$ and let $G=\gal(M/K)$. By Corollaries \ref{equality} and \ref{alltame}, the subgroups of $G$ defined in Proposition \ref{Hi} satisfy $H_i=\overline{H}_i$ for all $0\le i\le r$. Therefore,  $\left\{\delta_0,\dots,\delta_{r-1}\right\}$ is the underlying set of the multiset $\left\{v\left(\t-\sigma(\t)\right) \mid \sigma \in G, \ \sigma(\t)\ne\t\right\}$. It remains to find a concrete formula for the multiplicity $t_i$ of each value $\delta_i$.
	To this end, consider the chain of subgroups
\[
G=H_0\supset H_1\supset\cdots\supset H_{r-1}\supset H_{r}=\gal(M/K(\t)).
\]
The corresponding chain of fixed fields is $K=K(\al_0)\sub \cdots\sub K(\al_{r-1})\sub K(\t)=L$.

The natural action of $G$ on $\z(F)$ induces a bijection:
$$
G/H_r\lra \z(F),\qquad \sigma\longmapsto \sigma(\t).
$$
which   restricts to bijections:
$$
H_i/H_r\lra\z_i(F):=\left\{\t'\in \z(F)\mid v\left(\t-\t'\right)\ge\delta_i\right\}.
$$
Hence, the multiplicity $t_i=\#\z_i(F)-\#\z_{i+1}(F)$ is equal to:
\[
t_i=\#H_i/H_r-\#H_{i+1}/H_r=[L\colon K(\al_i)]-[L\colon K(\al_{i+1})]=(n/m_i)-(n/m_{i+1}).
\]
This ends the proof.
\end{proof}\e

 We may deduce from Theorem \ref{Valuesti} a nice relationship between weights and distances of the minimal polynomials  $\phi_i=\irr_K(\al_i)$, for all $0\le i< r$.  

\begin{corollary}\label{w-w}
Suppose that $\t$ is tame. Then,
\[
 d(\phi_i)-d(\phi_{i-1})=m_i\left(w(\phi_i)-w(\phi_{i-1})\right)\quad\mbox{for all }\ 0<i<r. 
\]
\end{corollary}

\begin{proof}
Let us fix some index $0<i<r$. By Lemma \ref{recurrence}, $[\al_0,\dots,\al_i]$ is an Okutsu sequence of $\al_i$. On the other hand, $\al_i$ is tame by Corollary \ref{alltame}. 
Theorem \ref{Valuesti}, applied to $\al_i$, yields the following equality  of multisets:
\[
\left\{v\left(\al-\al'\right) \mid \al'\in \op{Z}(\phi_i), \ \al'\ne\al\right\}=\left\{\delta_0^{t_{i,0}},\dots,\delta_{i-1}^{t_{i,i-1}}\right\},
\]
where \ $t_{i,j}=(m_i/m_j)-(m_i/m_{j+1})$ \ for all $0\le j< i$. Now, for every $\al'\in\z(\phi_i)$, $\al'\ne\al$, we have $v(\al-\al')\le\dta_{i-1}<\dta_i=v(\t-\al)$. Hence,
\[
v(\t-\al')=\min\{v(\t-\al),v(\al-\al')\}=v(\al-\al'). 
\]
 We deduce a formula for $w(\phi_i)$:
 \[
  m_iw(\phi_i)=\sum_{\al'\in\z(\phi_i)}v(\t-\al')=\dta_i+\sum_{\al'\ne\al}v(\al-\al')=t_{i,0}\dta_0+\cdots+t_{i,i-1}\dta_{i-1}+\dta_i.
 \]
For the weight of $\phi_{i-1}$ we have an analogous formula, so that
\[
  m_iw(\phi_{i-1})=\dfrac{m_i}{m_{i-1}}\,m_{i-1}w(\phi_{i-1})=\dfrac{m_i}{m_{i-1}}\left(t_{i-1,0}\dta_0+\cdots+t_{i-1,i-2}\dta_{i-2}+\dta_{i-1}\right).
\]
Straightforward computation shows that 
\[
 m_i\left(w(\phi_i)-w(\phi_{i-1})\right)=\dta_i-\dta_{i-1}.
\]
This ends the proof, because  $d(\phi_i)=\dta_i$ for all $i$, by the definition of an Okutsu sequence. 
\end{proof}\e

We end this section with an explicit formula for the main invariant $\delta(\t)$ in terms of the discrete invariants attached to a MLV chain of the valuation $v_F$.

By Theorems \ref{MLOk}, \ref{OkML} and $\ref{SF}$, the valuation $v_F$ admits an MLV chain 
$$
v\stackrel{\phi_0,\ga_0}\lra\  \mu_0\ \stackrel{\phi_1,\ga_1}\lra\  \mu_1\ \stackrel{\phi_2,\ga_2}\lra\ \cdots
\ \stackrel{\phi_{r-1},\ga_{r-1}}\lra\ \mu_{r-1} 
\ \stackrel{F,\infty}\lra\ \mu_{r}=v_F
$$
where $\ga_i=\mu_i(\phi_i)=v(\phi_i(\t))$ for all $i$. Since $F$ is defectless, all augmentations are ordinary.
Recall the secondary slopes introduced in Section \ref{subsecMLV}:
\[%\begin{equation}\label{gammai}
\lambda_0:=\ga_0, \qquad\lambda_i:=\ga_i-\mu_{i-1}(\phi_i),\quad 0<i\le r.
\]%\end{equation}

\begin{corollary}\label{lambdas2}
If $\t$ is tame, then \ $\delta_i=\lambda_0+\cdots+\lambda_i$,  for all \ $0\le i < r$.
\end{corollary}

\begin{proof}
For $i=0$, we have $\phi_0=x-\al_0$ and
$$\la_0=\ga_0= v(\phi_0(\t))=v(\t-\al_0)=\delta_0.$$
Hence, the claimed equality is equivalent to
\begin{equation}\label{desired}
\dta_i-\dta_{i-1}=\la_i=\ga_i-\mu_{i-1}(\phi_i).
\end{equation}
Now, $\ga_i=v_F(\phi_i)=m_iw(\phi_i)$. Also, since $\phi_{i-1}$ and $\phi_i$ are both key polynomials for $\mu_{i-1}$, Proposition \ref{weightMu} shows that
\[
 \mu_{i-1}(\phi_i)=\dfrac{m_i}{m_{i-1}}\,\mu_{i-1}(\phi_{i-1})=\dfrac{m_i}{m_{i-1}}\,\ga_{i-1}=\dfrac{m_i}{m_{i-1}}\,v_F\left(\phi_{i-1}\right)=m_i w(\phi_{i-1}).
\]
Thus, $\la_i=  m_i\left(w(\phi_i)-w(\phi_{i-1})\right)$ and the equality (\ref{desired}) follows from Corollary \ref{w-w}.
\end{proof}

\section{Appendix: some examples}\label{secExs}

In this Section, we display several examples illustrating the equivalence between Okutsu sequences of $\t$ and MLV chains of $v_F$. We address the following questions too:

\begin{itemize}
\item 	The computation of the relative invariants and the dependence of  these data on the choice of the generator $F$ of a fixed extensions $L/K$.
\item  The failure (or not) of the equality $\dta(\t)=\om(\t)$. 
\end{itemize}

In order to compute Krasner's constant, we shall use the technique of the \textbf{ramification polygon} \cite{RamP}.
If $L=K(\t)$ and  $F=\irr_K(\t)$, we may consider the polynomial 
\[
 \op{R}_F:=\dfrac1{\t^n}\,F(\t x+\t)\in L[x].
\]
If $F$ has roots $\t,\t_2,\dots,\t_n$, then $\op{R}_F$ has roots $(\t-\t_2)/\t,\dots (\t-\t_n)/\t$, besides the obvious root $0$. Therefore, the slopes of the Newton polygon of $\op{R}_F$ determine all values $v(\t-\t_i)$ and their multiplicities.

\subsection{Defectless examples}\label{subsecExDless}\mbox{\null}\e

\nn{\bf Example A1.} Let $K$ be the diadic field $\Q_2$, equipped with the $2$-adic valuation $v=\ord_2$. Consider the monic, irreducible polynomial
\[
 F=(x^2-2)^2+4x=x^4-2x^2+4x+4.
\]
Let $\z(F)=\{\t,\t_2,\t_3,\t_4\}$ and $L=K(\t)$.  We can obtain an MLV chain of $v_F$ as follows:
\[
 v\ \stackrel{x,\,1/2}\lra \ \mu_0\ \stackrel{\phi,\,5/4}\lra \ \mu_1\ \stackrel{F,\,\infty}\lra \ v_F,
\]
where $\phi=x^2-2$ and both augmentations are ordinary. Hence, $F$ is defectless and has depth two, with relative invariants:
\[
e_0=2,\ f_0=1; \quad  
e_1=2,\ f_1=1\quad\imp\quad e(L/K)=4,\ f(L/K)=1. 
\]
The corresponding Okutsu sequence is  
\[
[\al_0=0,\,\al_1=\sqrt2,\,\al_2=\t],\qquad  m_1=2. 
\]

Let us compute the main invariant. In the MLV chain we read that 
\[
\dta_0=v(\t)=v_F(x)=1/2,\qquad v_F(\phi)=5/4,\qquad \dta_2=v_F(F)=\infty.  
\]
At least one of the values  $v(\t+\sqrt2)$, $v(\t-\sqrt2)$ is smaller than  $3/2=v(2\sqrt2)$, because
\[
5/4=v(\t^2-2)=v(\t+\sqrt2)+v(\t-\sqrt2). 
\]
Since $\t+\sqrt2=\t-\sqrt2+2\sqrt2$, we deduce that 
\[
\dta(\t)=\dta_1= v(\t-\sqrt2)=v(\t+\sqrt2)=5/8.
\]
Since $e(L/K)=4$, we have $\sqrt2\not\in L$, so that $K(\al_1)\not\sub L$. Thus, the chain in (\ref{chainAlpha}) no longer holds beyond the tame case.

The ramification polygon of $F$ is one-sided of slope -1/6. Hence, for $i=2,3,4$ we have
\[
v\left((\t-\t_i)/\t\right)=1/6\ \mbox{ and }\ v(\t-\t_i)=\dfrac16+\dfrac12=\dfrac23.
\]
Therefore, $\om(\t)=2/3>\dta(\t)$.\bs

\nn{\bf Example A2. }In the same field $L$ of  Example  A1, let us choose an Eisenstein generator. For instance, the element \ $\eta:=(\t^2-2)/2$ \ satisfies \ $\eta^2=-\t$ \ and
\[
G:=\irr_K(\eta)=x^4-2x-2.
\] 
This new generator $G$ of the previous extension $L/K$ has depth one and $v_G$ admits the following MLV, consisting of  a single ordinary augmentation:
\[
v\ \stackrel{x,\,1/4}\lra \ \mu_0\ \stackrel{G,\,\infty}\lra \ v_G, \qquad e_0=4,\ f_0=1.
\]
The corresponding Okutsu sequence is 
\ $[\al_0=0,\,\alpha_1=\eta]$.

 Thus, the main invariant of $\eta$ is $\dta(\eta)=\dta_0=1/4$. The ramification polygon of  $G$ is one-sided of slope -1/12. Hence, for $i=2,3,4$ we have
\[
v\left((\eta-\eta_i)/\eta\right)=1/12\ \mbox{ and }\ v(\eta-\eta_i)=\dfrac1{12}+\dfrac14=\dfrac13,
\]
Therefore, $\om(\eta)=1/3>\dta(\eta)$.

Since $\al_0\in K$, this separable and defectless $\eta$ is a counterexample to \cite[Thm. 3.4]{OkS}.\bs

\nn{\bf Example B. }Let $k$ be an algebraically closed field of characteristic $p>0$. Let $t,y$ be indeterminates and consider $K=k(t)(\!(y)\!)$, equipped with the $y$-adic valuation $v=\ord_y$. Let $\mm$ be the maximal ideal of the valuation ring. Take
\[
F=x^p-y^{p-1}x-t\in \kx.
\]
Let $L=K(\t)$ for some $\t\in\z(F)$. This polynomial $F$ is irreducible and $f(L/K)=p$, because $F\equiv x^p-t\md{\mm}$  is irreducible modulo $\mm$. In particular, $Lv/K\!v$ is purely inseparable.

A MLV chain of $v_F$ and the corresponding Okutsu sequence would be:
\[
v\ \stackrel{x,\,0}\lra \ \mu_0\ \stackrel{F,\,\infty}\lra \ v_F, \qquad [\al_0=0,\al_1=\t],
\]
where the unique augmentation is ordinary.
The main invariant of $\t$ is $\dta(\t)=\dta_0=0$. 

Since $\z(F)=\t+y\F_p$, we see that Krasner's constant of $\t$ is $\om(\t)=1>\dta(\t)$.

On  the other hand, this example is another counterexample to \cite[Thm. 3.4]{OkS}.

\subsection{Defect examples}\label{subsecExDefect}
Let $k$ be an algebraically closed field of odd characteristic $p$.
For an indeterminate $y$, consider the field $k(\!(y^\Q)\!)$ of power series in $y$ with rational exponents. 
The \textbf{support} of a power series is the following set: 
$$s=\sum\nolimits_{q\in\Q}a_qy^q\ \imp\ \supp\left(s\right)=\{q\in\Q\mid a_q\ne0\}\sub\Q.$$
The Hahn field $\Ha\subset k(\!(y^\Q)\!)$ consists of all power series with well-ordered support. It is an algebraically closed field. The valuation $v=\ord_y$ on $k(y)$ admits a unique extension to $\Ha$, defined by \ $v(s):=\min(\supp(s))$. 
Clearly, 
\[
v\Ha=\Q \quad\mbox{ and }\quad \Ha v=k. 
\]

Let $K$ be be the perfect hull of $k(\!(y)\!)$ in $\Ha$. That is,
\[
K=\bigcup\nolimits_{n\in\N} k(\!(y^{1/p^n})\!).
\]
Let $\kb\sub \Ha$ be the algebraic closure of $K$ in $\Ha$.
Any finite subextension $K\sub L\sub\kb$ with $[L\colon K]=p^m$ is necessarily immediate, with $d(L/K)=p^m$.\bs

\nn{\bf Example C. }Consider the Artin-Schreier polynomial
\[
F=x^p-x-y^{-1}\in\kx.
\]
We can exhibit a concrete root $\al$ of $F$ in $\Ha$, namely
\[
\al=\sum_{i=1}^\infty y^{-1/p^i}\in\Ha,\qquad v(\al)=-1/p.
\]
Clearly,  $\al\not\in K$ because the denominators of the exponents of $y$ are unbounded. Since $\z(F)=\al+\F_p$, we see that $F$ is irreducible in $\kx$ and Krasner's constant of $\al$ is $\om(\al)=0$. 

Consider the partial sums of the series defining $\al$:
\[
\al_m:=\sum_{i=1}^m y^{-1/p^i}\in K,\quad\mbox{for all }\ m\in\N.
\] 
An Okutsu sequence of $\al$ could be, for instance:
\[
\left[A_0=\left\{\al_m\right\}_{m\in\N},\,A_1=\{\al\}\right].
\] 
Thus, $\al$ has depth one and a MLV chain of $v_F$ would consist of a single limit augmentation. The main invariant of $\al$ is
\[
\dta(\al)=\sup\{v(\al-\al_m)\mid m\in\N)\}=0=\om(\al).
\]\e

\nn{\bf Example D1. }Let $L=K(\al,y^{1/2})$, which has degree $2p$ over $K$, with $e(L/K)=2$ and $d(L/K)=p$.
Take $\t=y^{1/2}\al$ as a generator of $L$. The minimal polynomial of $\t$ over $K$ is:
\[
F=\left(x^p-y^{(p-1)/2}x\right)^2-y^{p-2}.
\]
An Okutsu sequence of $\t$ is
\[
\left[A_0=\{0\},\,A_1=\left\{y^{1/2}\al_m\right\}_{m\in\N},\,A_2=\{\t\}\right], \qquad m_1=2.
\] 
Thus, $F$ has depth two and the MLV chains of $v_F$ will have an  ordinary augmentation, followed by a limit augmentation. 
The main invariant of $\t$ is:
\[
\dta(\t)=\sup\left\{v\left(\t-y^{1/2}\al_m\right)\mid m\in\N\right\}=1/2.
\] 

It is easy to check that the roots of $F$ are:
\[
\z(F)=y^{1/2}\left(\al+\F_p\right)\cup -y^{1/2}\left(\al+\F_p\right).
\]
Hence, $\om(\t)=1/2=\dta(\t)$.\bs

\nn{\bf Example D2. }For the same field $L=K(\al,y^{1/2})$ of the previous example, take  $\eta=\al+y^{1/2}$ as a generator over $K$. The minimal polynomial of $\eta$ over $K$ is:
\[
G=\left(x^p-x-y^{-1}\right)^2-y\left(y^{(p-1)/2}-1\right)^2.
\]
An Okutsu sequence of $\t$ is
\[
\left[A_0=\left\{\al_m\right\}_{m\in\N},\,A_1=\{\al\},\,A_2=\{\t\}\right], \qquad m_1=p.
\] 
Thus, $G$ has depth two and the MLV chains of $v_G$ will have a limit augmentation, followed by an  ordinary augmentation. 
The main invariant of $\eta$ is: $\dta(\eta)=v\left(\eta-\al\right)=1/2$.

The roots of $G$ are: $\z(G)=\al+\F_p\pm y^{1/2}$, so that
 $\om(\eta)=1/2=\dta(\eta)$.\bs
 
 \nn{\bf Example E. }For the same $\al\in\z(x^p-x-y^{-1})$ of the previous examples, take 
 \[
  \t:=\sum_{i=1}^\infty\left(y^{-1}\al\right)^{1/p^i}\in\Ha.
 \]
Since $\t$ satisfies $\t^p-\t=y^{-1}\al$, its minimal polynomial over $K$ is
\[
F=\left(x^p-x\right)^p-y^{1-p}\left(x^p-x\right)-y^{-1-p}
\]
The extension $L=K(\t)$ has $[L\colon K]=d(L/K)=p^2$. An Okutsu sequence of $\t$ is
\[
\left[A_0=\left\{(y^{-1}\al_m)^{1/p}\right\}_{m\in\N},\,A_1=\left\{\sum_{i=1}^m\left(y^{-1}\al\right)^{1/p^i}\right\}_{m\in\N},\,A_2=\{\t\}\right], \qquad m_1=p.
\] 
Thus, $F$ has depth two and the MLV chains of $v_F$ will have two consecutive limit augmentations. 
The $\dta$-invariants of the sequence are:
\[
 \dta_0=-1/p,\qquad \dta_1=\dta(\t)=0. 
\]
Since $\z(F)=\t+\F_p+\F_p\al$, we have $\om(\t)=0$ as well.

\end{document}